\documentclass[a4paper,11pt]{article}
\usepackage{bbold}
\usepackage{mathrsfs}
\usepackage{mathrsfs}
\usepackage{amsmath}
\usepackage{amsfonts}
\usepackage{amsthm}
\usepackage{times}
\usepackage{hyperref}
\usepackage[numbers,sort&compress]{natbib}

\usepackage{amssymb, graphicx, epsfig}

\usepackage{latexsym}

\hypersetup{
    bookmarks=true,         % show bookmarks bar?
    unicode=false,          % non-Latin characters in Acrobat’s bookmarks
    pdftoolbar=true,        % show Acrobat’s toolbar?
    pdfmenubar=true,        % show Acrobat’s menu?
    pdffitwindow=true,      % page fit to window when opened
    pdftitle={My title},    % title
    pdfauthor={Author},     % author
    pdfsubject={Subject},   % subject of the document
    pdfnewwindow=true,      % links in new window
    pdfkeywords={keywords}, % list of keywords
    colorlinks=true,       % false: boxed links; true: colored links
    linkcolor=blue,          % color of internal links
    citecolor=blue,        % color of links to bibliography
    filecolor=blue,      % color of file links
    urlcolor=blue           % color of external links
}

\newtheorem{con}{Conjecture}[section]
\theoremstyle{definition}
\newtheorem{rem}{Remark}[section]

\author{Johan~Bj\"orklund\footnote{Department of Mathematics, Uppsala University, SE},~~Cecilia~Holmgren\footnote{Department of Pure Mathematics and Mathematical Statistics, University of Cambridge, UK} \footnote{corresponding author: Cecilia Holmgren, DPMMS, Centre for Mathematical Sciences, Wilberforce Road, Cambridge CB3 0WA, UK, C.Holmgren@dpmms.cam.ac.uk}}
\title{Counterexamples to a Monotonicity Conjecture for the Threshold Pebbling Number}

\begin{document}\maketitle
\begin{abstract}
Graph pebbling considers the problem of transforming configurations of discrete pebbles to certain target configurations on the vertices of a graph, using the so-called \emph{pebbling move}.
This paper provides counterexamples to a monotonicity conjecture stated by Hurlbert et al. \cite{Hurlbert1, Hurlbert2, Hurlbert3, Hurlbert4} concerning the \emph{pebbling number} compared to the \emph{pebbling threshold}.
\end{abstract}
\textbf{Keywords:}combinatorics, probability theory, graph theory, graph pebbling, pebbling number, pebbling threshold
\section{Introduction}
Graph pebbling considers the problem of transforming configurations of discrete pebbles to certain target configurations on the vertices of a graph, using the so-called pebbling move.

Historically graph pebbling was first suggested by Lagarias and Saks in attempt to answer a number-theoretic question by Erd\H os and Lemke concerning zero-sum sequences of elements from a finite group, see \cite{Hurlbert2,Hurlbert3}. However, the concept of graph pebbling was first introduced in the literature by Chung \cite{Chung} who defined the  pebbling number $\pi(G)$ of a connected graph $G$. Today the area is very active,  with many open problems and conjectures.

Czygrinow, Eaton, Hurlbert and Kayll \cite{Hurlbert1} introduced a probabilistic pebbling model  and defined the so-called pebbling threshold $\tau(G)$. Our aim in this paper is to provide counterexamples to a monotonicity conjecture stated by Hurlbert et al. in \cite{Hurlbert1,Hurlbert2,Hurlbert3, Hurlbert4} relating the pebbling numbers $\pi(G_n)$ to the pebbling thresholds $\tau(G_n)$ for graph sequences.

A \emph{pebbling move} on a graph consists of removing two pebbles from one vertex and placing one on an adjacent vertex  (the second removed pebble is discarded from play).
A configuration of $t$ pebbles on the vertices of a graph $G$ is \emph{solvable} if for any vertex $v$ of $G$, it is possible after a series of pebbling moves to reach a new configuration so that $v$ has one or more pebbles. A configuration which is not solvable is said to be \emph{unsolvable}.
 The \emph{pebbling number} $\pi(G)$ is the smallest $t$ such that all initial configurations of $t$ pebbles on the vertices of the graph $G$ is solvable.

In the probabilistic pebbling model introduced by Czygrinow et al. \cite{Hurlbert1} the pebbling configuration is selected uniformly at random from the set of all possible configurations with $t$ pebbles. 
(This is one of many possible random models, e.g., one could consider a model where each of the pebbles is uniformly at random placed on a vertex of $G$.) 
%Heuristically the \emph{pebbling threshold} is the smallest number so that if one takes $t$ much larger than the threshold a random configuration is a.s. solvable, and if we take $t$ much smaller than the threshold a random configuration is a.s. unsolvable.

Below  we define the \emph{pebbling threshold} on a \emph{graph sequence}, see \cite{Brightwell}.  
Consider a graph sequence $\mathcal{G}=(G_1,G_2,\dots, G_n,\dots)$, where $G_n$ has vertex set $[n]$.
Let $\mathbb{N}$ be the set of non-negative integers and let $C_n:[n]\rightarrow \mathbb{N}$ denote a configuration of pebbles on $n$ vertices. For a function $t=t(n)$ we let 
$\mathcal {D}(G_n,t)$ be the probability space of all configurations $C_n$ of size $t=\sum_{i\in [n]}C_n(i)$ in $G_n$,  with each configuration having the same probability i.e., $1/\binom{n+t-1}{t}$.
We write $P(G_n;t)$ the probability that a configuration chosen uniformly at random from $\mathcal {D}(G_n,t)$ is solvable in $G_n$. For $\alpha\in (0,1)$ we define  \[\tau_\alpha(n)=\tau_{\alpha}(G_n):=\min\{t:~~P(G_n;t)\geq \alpha\}\] and call it a threshold function for $G_n$. As is customary we consider $\alpha=\frac{1}{2}$ and write for simplicity 
$\tau:=\tau_{\frac{1}{2}}$. 
Bekmetjev, Brightwell, Czygrinow and Hurlbert \cite{Brightwell} showed that for any sequence  $\omega=\omega(n)$ tending to infinity we have
\[P(G_n,\tau\omega)\rightarrow 1~~ \mathrm{and}~~P(G_n,\tau/\omega)\rightarrow 0,~~~~\mathrm{as}~~n\rightarrow\infty.\]

 %$t_{\frac{1}{2}}\in \tau(\mathcal{G})$ implying that there exist threshold functions for all graph sequences.
%When we write $\tau(\mathcal{G})<\tau(\mathcal{H})$ for graph sequences $\mathcal{G}=(G_1,G_2,\dots,G_n,\dots)$ and $\mathcal{H}=(H_1,H_2,\dots,H_n,\dots)$, we mean that for all threshold functions $t(n)\in \tau(\mathcal{G})$ and  $\widehat{t}(n)\in \tau(\mathcal{H})$ it holds that $\lim_{n\rightarrow\infty}\frac{t(n)}{\widehat{t}(n))}=0$. Furthermore, when we write $\tau(\mathcal{G})=\tau(\mathcal{H})$ we mean that there exists a constant $C>0$ such that that for all threshold functions $t(n)\in \tau(\mathcal{G})$ and  $\widehat{t}(n)\in \tau(\mathcal{H})$ it holds that $0<\lim_{n\rightarrow\infty}\frac{t(n)}{\widehat{t}(n)}\leq C$.

In \cite[Question 2.3]{Hurlbert1}, \cite[Conjecture 4.5]{Hurlbert2}, \cite[Conjecture 8.4]{Hurlbert3} and  \cite[Section 4.2 p.20]{Hurlbert4} Hurlbert et al. stated the following conjecture concerning the relationship of the pebbling number to the pebbling threshold.
\begin{con}{\emph{Hurlbert et al.}}\label{conjecture} If $\mathcal{G}=(G_1,G_2,\dots,G_n,\dots)$ and $\mathcal{H}=(H_1,H_2,\dots,H_n,\dots)$ are graph sequences such that $\pi(G_n)\leq \pi(H_n)$, then for $\alpha\in (0,1)$ it holds that $\tau_{\alpha}(G_n)\in O(\tau_{\alpha}(H_n))$.
\end{con}
In Section \ref{2} we disprove this conjecture by constructing two counterexamples.

\section{Disproving Conjecture \ref{conjecture}}\label{2}

\begin{figure}[h]
\begin{center}
\includegraphics[scale=0.20]{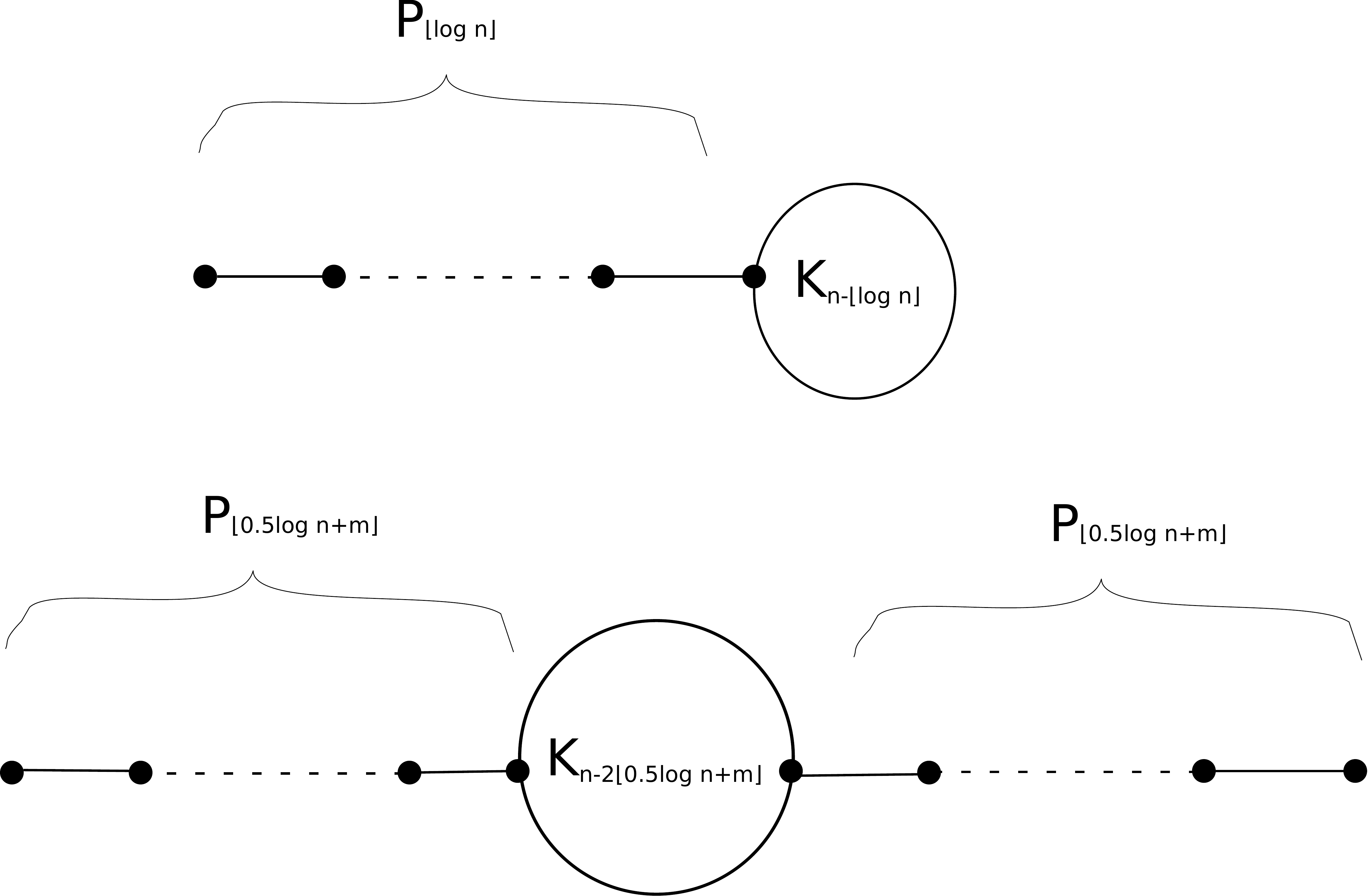}
\end{center}
\caption{\small{\textit{For the graph sequences $\mathcal{G}=(G_1,G_2,\dots,G_n,\dots)$ and $\mathcal{H}=(H_1,H_2,\dots,H_n,\dots)$ it holds that $\pi(G_n)<\pi(H_n)$ but $\tau(G_n)\not\in O(\tau(H_n))$.}}}
\label{gnhn}
\end{figure}

Let $P_n$ denote the \emph{path} with $n$ vertices and let $K_n$ denote the \emph{complete graph} with $n$ vertices. When we write $\log$ we mean $\log_2$. Let $G_n$ be a graph consisting of a path with $\lfloor \log n\rfloor$ vertices such that the last vertex  is connected by one edge to a complete graph with 
$n-\lfloor \log n\rfloor$ vertices. 
 Let further $H_n$ be a graph consisting of two paths with $\lfloor \frac{1}{2}\log n+m\rfloor$ vertices (where we can choose $m$ for example equal to 1000) such that both paths have their endpoints connected by one edge to a complete graph with $n-2(\lfloor \frac{1}{2}\log n+m\rfloor)$ vertices. See Figure \ref{gnhn} where $G_n$ and $H_n$ are illustrated.

We will now see that $H_n$ has a larger pebbling number than $G_n$.
We see that $G_n$ and $H_n$ could be regarded as paths with roughly $\log n$ and $\log n+m$ vertices, respectively, where one vertex is replaced by an entire complete graph. The vertices that are hardest to reach are the vertices at the ends of the paths. 
We claim that $3n\geq \pi(G_n)$.
Indeed, in order to move pebbles to the complete graph it is enough to have 2 pebbles in the vertex which joins the complete graph to the path and thus $2^{\log n}=n$ pebbles is enough, starting in the other endpoint of the path. Also, if we want to move pebbles from the complete graph to the other endpoint of the path a worst configuration is to place 3 pebbles in each vertex in the complete graph. Since we need at most $2^{\log n}$ pebbles in the path for the path to be pebbled, we see that $3n$ is indeed an upper bound for the pebbling number of $G_n$, as claimed.
On the other hand, for $H_n$ it is obvious that we need a pebbling number which is at least as large as the pebbling number for the path with roughly $\log n+m$ vertices. Thus, it is obvious choosing $m=1000$ that $\pi(H_n)\geq 100n$.

Now we will see that the graph sequence $\mathcal G=(G_1,G_2,\dots,G_n,\dots)$ has a larger threshold than the graph sequence $\mathcal H=(H_1,H_2,\dots,H_n,\dots)$. 
We first explain why it is at least intuitively reasonable that $\tau (G_n)\not\in O(\tau(H_n))$.
In the path $P_n$, the vertices which are the hardest to pebble are the endpoints, and the maximal number of pebbles that are needed to pebble such a vertex is $2^{n-1}$, when all pebbles are in the opposite endpoint in the initial configuration. This is the scenario we are trying to achieve by our construction of $G_n$, forcing the starting configuration to (with high probability) have almost all pebbles in the complete graph (containing almost all the vertices of $G_n$), so that the opposite endpoint of the path will be hard to pebble.
 
 On the other hand, by our construction of $H_n$,  where the (large) complete graph is now put in the middle of the path the initial configuration is forced to (with high probability) have almost all pebbles in the middle of the path.  We note that if all pebbles are placed in the vertex which is closest to the middle of $P_n$ then  $2^{\lfloor\frac {n}{2}\rfloor}$ pebbles are enough to pebble both endpoints.
Thus, we claim that the threshold for $\mathcal G$ should be higher than the threshold for $\mathcal H$.
We will now formally prove that this claim holds. 

We first recall some standard asymptotic notation.  We write $f\in O(g)$ 
(equivalently $g\in \Omega (f)$) when there are positive constants $c$ and $k$ such that $f(n)/g(n)<c$ for all $n>k$, and we denote $\Theta(g)$ for $O(g)\cap \Omega (g)$.

Returning to the sequence $\mathcal G=(G_1,G_2,\dots,G_n,\dots)$, suppose that we add at most $n^{0.99}$ pebbles to $G_n$. Then by symmetry, the expected number of pebbles that are distributed to the path $P_{\lfloor\log n\rfloor}$ is $O(\frac{n^{0.99}\log n}{n})$. Hence, the Markov inequality implies that the probability that there are any pebbles in the path tends to 0. Note that even if all $n^{0.99}$ pebbles in the complete sub-graph were to lie in the same vertex we would not be able to reach the last vertex in the path.
Thus,  $\tau(G_n)\in\Omega (n^{0.99})$.
%Since this is less than the path threshold $O(\lg n 2^{c\log\log n})$ the path can not pebble itself, meaning that it need to move pebbles from the complete graph in the end of the path in $G_n$.
%Choosing $C$ large enough we see that we can pebble, that is we can move $2^{\lg_2 n-1}$ pebbles from one endpoint to vertices close to the other end point which is the worst case scenario. 
%Now we will see that the threshold for $G_n$ is $\Theta(n)$.

%We will show that if we take less than $n/1000$ pebbles randomly we will not pebble $G_n$ with probability tending to 1.

%in the complete graph we can not move the pebbles to the opposite end point of the path. 

%There has to be an empty vertex with high probability in a vertex which is quite close to the opposite end point say an interval which is 1 procent of the total length of the path.
%If we 
Now we return to the sequence $\mathcal H=(H_1,H_2,\dots,H_n,\dots)$. We first observe that for $H_n$, since there are two paths of length $\lfloor\frac{\log n}{2}+m\rfloor$ instead of one longer path as in $G_n$, it is enough to move $b\sqrt n$ pebbles (where $b$ is a constant) to one vertex in the complete sub-graph to pebble every vertex of $H_n$. 
Our next aim is to show that $\tau(H_n)\in O(n^{0.8})$. 
Given a configuration $C$ of pebbles in a graph $G$, we define the \emph{number of birthdays} $B(C)$ as the number of pebbles 
in the configuration $C^{'}$ obtained by removing one pebble from each vertex of $G$ assigned at least one pebble by $C$.
 after having removed one pebble from each vertex in $G$ containing a pebble. %For $t=n^{0.8}$ there are at most $j\binom n t$ configurations such that $B(C)\leq j$ when $j\leq n^{0.55}$. 
We first note that the number of configurations such that $B(C)=k$ is $\binom {n}{t-k}\binom{t-k+k-1}{k}$ by first choosing $t-k$ vertices where we place one pebble, and then distribute $k$ pebbles to these $t-k$ vertices. Let $a(k):=\binom {n}{t-k}\binom{t-1}{k}$. We write $S(j):=\sum_{k=0}^{j}a(k)$. Note that $S(t-1)=\binom{n+t-1}{t}$, since each possible configuration with $t$ pebbles is calculated exactly once in the sum $S(t-1)$.  If $B(C)\geq 2c\sqrt {n}$ we can move at least $c\sqrt{n}$ distinct pebbles. Hence, it is enough to prove that for $t=n^{0.8}$ we have \begin{align}\label{0}\dfrac{S(n^{0.55})}{S(t-1)}\rightarrow 0,~~~~\mathrm{as}~~n\rightarrow \infty,
\end{align}since this would imply that we get at least $n^{0.55}$ birthdays with probability tending to 1.
To prove that (\ref{0}) holds we first observe that \begin{align}\label{boundy}\frac{a(k-1)}{a(k)}=\frac{k(n-t+k)}{(t-k)(t-k+1)}<1,\end{align} when $nk<t^2-kt+t-k$.
Hence, for $t=n^{0.8}$ we have
\begin{align}\label{1}\frac{S(n^{0.55})}{S(t-1)}\leq\frac{S(n^{0.55})}{S(n^{0.6})-S(n^{0.55})}\leq \frac{n^{0.55}a(n^{0.55})}{a(n^{0.55})(n^{0.6}-n^{0.55})}\in O(n^{-0.05}),
\end{align}
proving the limit in (\ref{0}).

So for $\mathcal G$ we have a threshold which is at least $n^{0.99}$ and for $\mathcal H$ we have a  threshold which is at most $n^{0.8}$. 
Thus, $\tau(G_n)\not\in O(\tau(H_n))$  but $\pi(G_n)<\pi(H_n)$ disproving the conjecture.

\begin{rem}\label{remark1}In relation to Conjecture \ref{conjecture}, Czygrinow et al. \cite{Hurlbert1} suggested that it is possible that the conjecture only holds with the additional hypothesis that $\pi(G_n)$ is significantly smaller than $\pi(H_n)$. They suggested that it might be possible that one needs to add the condition that  $\pi(H_n)-\pi(G_n)\rightarrow \infty$ or $\lim \sup_{n\rightarrow \infty}\pi(G_n)/\pi(H_n)<1$ for the conjecture to hold.
With our counterexample, now choosing $m=\log\log n$ instead of 1000, we can easily see that even if $\lim_{n\rightarrow \infty}\pi(G_n)/\pi(H_n)=0$ (which obviously now is the case) it still holds that $\tau(G_n)\not\in O(\tau(H_n))$ . Thus, our 
counterexample also disproves the conjecture even with this additional hypothesis.
\end{rem}

\begin{rem}

%This remark is related to Remark \ref{remark1}. 
%In \cite[Theorem 1.1]{Brightwell} it was shown that for any  graph sequence and for any $\epsilon >0$ the threshold is $\Omega{n^{1/2}}\cap o(n^{1+\epsilon})$.
It is obvious that the lowest threshold for graph sequences is $\Omega({n^{1/2}})$, the threshold for complete graphs, and the highest is the threshold for paths, i.e., $O(n 2^{c\sqrt{\log n}})$ for any constant $c>1$ (which is in particular $o(n^{1+\epsilon})$ for all $\epsilon >0$), see \cite{Godbole}.

\begin{figure}[h]
\begin{center}
\includegraphics[scale=0.18]{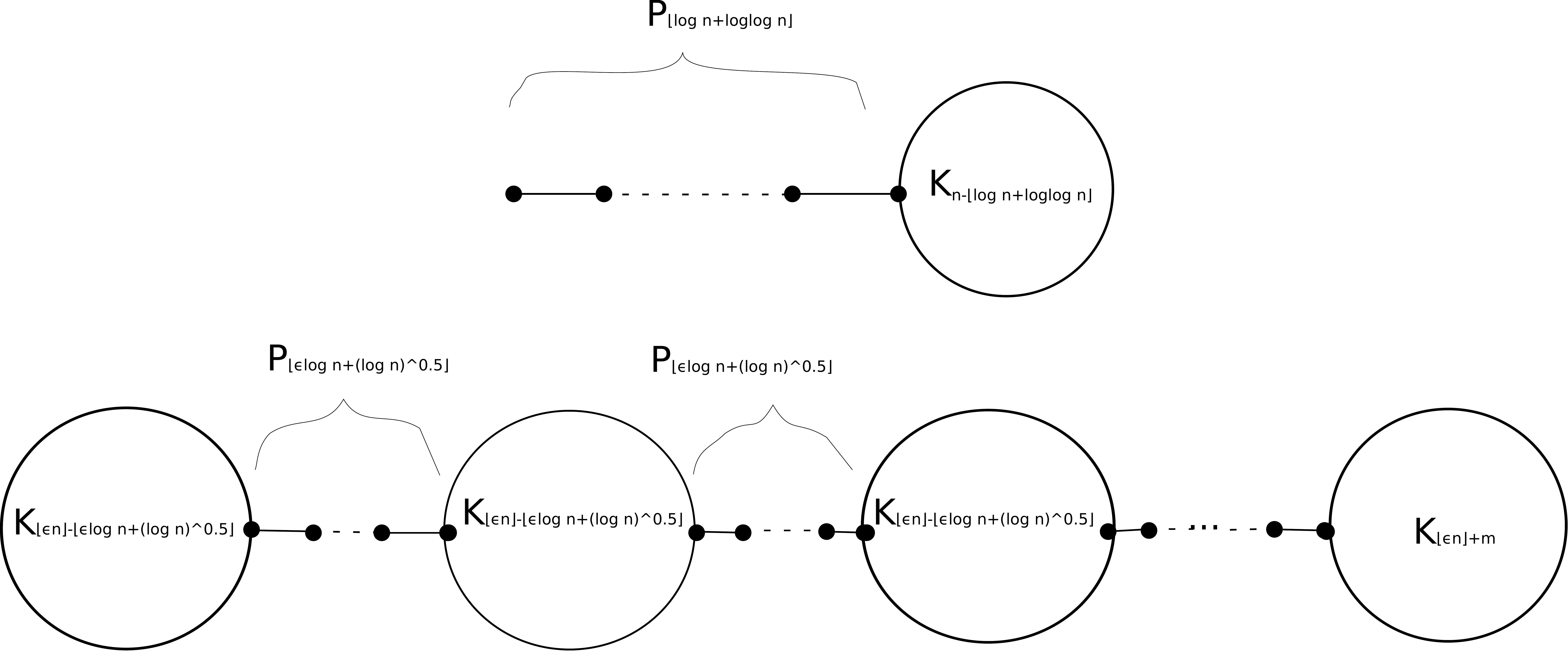}
\end{center}
\caption{\small{\textit{For the graph sequences $\mathcal{G}=(G_1,G_2,\dots,G_n,\dots)$ and $\mathcal{H}=(H_1,H_2,\dots,H_n,\dots)$ we have $\pi(G_n)<\pi(H_n)$. However, it holds that $\tau(G_n)\in\Omega(n\log n)$ and $\tau(H_n)\in O(n^{\frac {1}{2}+\epsilon})$ for any $\epsilon>0$, i.e., $\tau(G_n)$ and $\tau(H_n)$ are close to the largest respectively smallest threshold of a graph sequence.}}}
\label{gnhn2}
\end{figure}

Consider a modification of our counterexample such that $\mathcal{G}=(G_1,G_2,\dots,G_n,\dots)$ is the same sequence as in our example above except from that we let the path have $\lfloor \log n+\log\log n\rfloor $ vertices and the complete graph have $n-\lfloor \log n+\log\log n\rfloor $ vertices. 
It is easy to see from our arguments above that the order of the threshold  $\tau(G_n)$ is larger than $n$.
For the graph sequence $\mathcal{H}=(H_1,H_2,\dots,H_n,\dots)$, we modify our example by first considering a number of complete graphs $S_i$ indexed by $i$, each  with $\lfloor \epsilon n\rfloor-\lfloor \epsilon \log n+\sqrt {\log n}\rfloor$ vertices. In analogy with our earlier construction, we attach a path containing $\lfloor \epsilon \log n+\sqrt {\log n}\rfloor$ vertices with one endpoint in $S_i$ and one endpoint in $S_{i+1}$ for all $i$.
We have roughly $\frac {1}{\epsilon}$ sub-graphs consisting of a path and a complete graph connected to it. 
(When we don't have enough vertices to construct a new sub graph consisting of a path and a complete graph, we let the last complete graph get the rest of the vertices.) 
For an illustration of these modified graph sequences, see Figure \ref{gnhn2}.

Obviously $\lim_{n\rightarrow \infty}\pi(G_n)/\pi(H_n)=0$.
However, the threshold $\tau(H_n)$ is at most $n^{\frac{1}{2}+\epsilon}$ by  calculations analogous to those in (\ref{0}), (\ref{boundy}) and (\ref{1}), changing $n^{0.55}$ to $n^{\epsilon}$ (corresponding to the number of birthdays), $n^{0.6}$ to $n^{2\epsilon}$ and $t=n^{0.8}$ to $t=n^{\frac{1}{2}+\epsilon}$.
 Thus, this modification of our counterexample shows that we can find graph sequences such that $\lim \sup_{n\rightarrow \infty}\pi(G_n)/\pi(H_n)=0$, but where $\tau(G_n)$ is close to the largest threshold of a graph sequence respectively $\tau(H_n)$ is close to the smallest threshold of a graph sequence.

\end{rem}

\subsection*{Acknowledgement:}
We thank Professor B\'ela Bollob\'as for introducing us to the problem area and for his helpful comments on the manuscript.

\end{document}